\documentclass[a4paper,10pt,reqno]{amsart}
\usepackage{amssymb} \usepackage{amsfonts}
\usepackage{amsmath,amsthm}
\usepackage[T1]{fontenc}

\usepackage{enumitem}
\setlist[enumerate,1]{label=(\roman*),leftmargin=2em}

\newtheorem*{mtheo}{Main Theorem}
\newtheorem{theo}{Theorem}[section]
\newtheorem{cor}{Corollary}[section]
\newtheorem{prop}{Proposition}[section]
\theoremstyle{definition}
\newtheorem{defi}{Definition}[section]
\newtheorem{ex}{Example}[section]
\theoremstyle{remark}

\DeclareMathOperator{\Z}{\mathbb{Z}}
\DeclareMathOperator{\R}{\mathbb{R}}
\DeclareMathOperator{\F}{\mathbb{F}}
\DeclareMathOperator{\Id}{I}
\DeclareMathOperator{\GL}{GL}

\usepackage{marginnote}
\usepackage{color}
\usepackage{changepage}
\newcounter{rlcheck}
\strictpagecheck
\usepackage[normalem]{ulem}
\makeatletter
\def\redsout{%
	\ULdepth=.5ex%
	\ifmmode \ULdepth=1.0ex \fi%
	\bgroup\markoverwith{\textcolor{red}{\rule[\ULdepth]{2pt}{1.5pt}}}\ULon%
}%
\makeatother
\begin{document}

\title[Spin-structures on real Bott manifolds with K\"{a}hler structure]{Spin-structures on real Bott manifolds with K\"{a}hler structure}

\author{Anna G\k{a}sior \and Rafa{\l} Lutowski}

\address{Anna G\k{a}sior:
Maria Curie-Sk{\l}odowska University,
Institute of Mathematics,
pl. Marii Curie-Sk{\l}odowskiej 1,
20-031 Lublin, Poland}
\email{anna.gasior@umcs.pl}

\address{Rafa{\l} Lutowski:
    Institute of Mathematics,
    University of Gda\'nsk,
    ul. Wita Stwosza 57,
    80-308 Gda\'nsk, Poland}
\email{rafal.lutowski@ug.edu.pl}

\begin{abstract}
In 2011 Ishida gave a criterion for existence of K\"{a}hler structures on real Bott manifolds. For those manifolds we present necessary and sufficients condition for admitting spin structures.
\end{abstract}

\maketitle

\section{Introduction }

Let $\Gamma$ be a fundamental group of an $n$-dimensional \emph{real Bott manifold} $M$. From \cite{KM09} we know that $\Gamma$ is a \emph{semi-diagonal} Bieberbach group (and $M$ is a \emph{semi-diagonal} flat manifold). This means that:
\begin{enumerate}
\item $\Gamma$ fits into a short exact sequence
\begin{equation}
\label{wzor21}
0\longrightarrow \Z^n\stackrel{\iota}{\longrightarrow} \Gamma \stackrel{\pi}{\longrightarrow} C_2^d \longrightarrow 1,
\end{equation}
for some $1 \leq d \leq n$, where $C_2$ is the cyclic group of order $2$;
\item the image of $\varrho \colon C_2^d\to \GL(n,\Z)$ is a group of diagonal matrices, where
$$\varrho_g(z)=\iota^{-1}(\gamma\iota(z)\gamma^{-1}),$$
for $z \in \Z^n, g\in C_2^n$ and $\gamma\in\Gamma$ is such that $\pi(\gamma)=g$.
\end{enumerate}
Note that in the particular case of real Bott manifolds we can take $d=n$.

Up to diffeomorphism, $M$ is determined by a certain strictly upper triangular matrix $A$ with coefficients in $\F_2$.
We call $A$ a \emph{Bott matrix} and we denote the manifold $M$  by $M(A)$. In \cite[Theorem 3.1]{I11} Ishida gave a necessary and sufficient condition for existence of the K\"{a}hler structure on $M(A)$: $A$ is of even dimension, say $2n$, and one can split its columns into $n$ pairs of equal ones.

Let $M$ be a real Bott and K\"ahler manifold, an \emph{RBK-manifold} for short. In this note we examine the existence of spin structures on $M$. We would like to mention that the general condition for existence of spin structures on semi-diagonal flat manifolds is considered in \cite{GS14} and \cite{G17}. However, in the specific case of RBK-manifolds, this condition can be formulated in much simpler -- purely combinatorial -- form. Namely, if $\tilde A$ is obtained from $A$ by a removal of one column from each pair of equal ones and $\tilde{S}_i$ is the sum of elements from the $i$-th row of $\tilde A$, then we get:

\begin{mtheo}
An RBK-manifold $M(A)$, of dimension $2n$, admits a spin structure if and only if for every $1 \leq i \leq 2n$
$$\tilde{S}_i=1\Longrightarrow A^{(i)}=0,$$
where $A^{(i)}$ denotes the $i$-th column of the matrix $A$.
\end{mtheo}

We prove the above theorem in Section \ref{sec:main}. We take advantage of the description of spin semi-diagonal flat manifolds presented in \cite{LPPS19}. For the convenience of readers it is recalled in Section \ref{sec:diagonal}.

\section{New definition of real Bott manifold}
\label{sec:diagonal}

\indent In this section we recall methods introduced in \cite{PS16} and developed in \cite{LPPS19}.
Let $S^1$ be a unit circle in $\mathbb C$ and we consider automorphisms $g_i:S^1\to S^1$ given by
\begin{equation}
g_0(z)=z,\;\; g_1(z)=-z,\;\; g_2(z)=\bar z,\;\; g_3(z)=-\bar z,
\end{equation}
for all $z\in S^1$. If we identify $S^1$ with $\R/ \Z$, then for each ${[t]\in \R/ \Z}$  we get
\begin{equation}\label{wzor23}
g_0([t])=[t],\;\; g_1([t])=\left[t+\frac12\right],\;\; g_2([t])=[-t],\;\; g_3([t])=\left[-t+\frac12\right].
\end{equation}
Let ${D}=\langle g_i:i=0,1,2,3\rangle$. Then ${D}\cong{C_2\times C_2}$ and $g_3=g_1g_2$. We define an action ${D}^n$ on $T^n$ by
\begin{equation}
(t_1,\ldots,t_n)(z_1,\ldots,z_n)=(t_1z_1,\ldots,t_nz_n)
\end{equation}
for $(t_1,\ldots,t_n)\in{D}^n$ and $(z_1,\ldots,z_n)\in T^n=\underbrace{S^1\times\ldots\times S^1}_{n}$.

By taking any $d$ generators of the group $C_2^d\subseteq D^n$, we define a $(d\times n)$-matrix with entries in $D$. This in turn defines a matrix with entries in the set $P=\{0,1,2,3\}$, under the identification $i\leftrightarrow g_i$ for $i=0,1,2,3$. We call it a $P$-matrix of $C_2^d$. Note that under the above identification, $P$ has the natural structure of a vector space over $\F_2$.

Although a group $C_2^d\subseteq D^n$ can have many $P$-matrices in general, every such a matrix $E$ encodes some important properties of its action the torus $T^n$. Namely, $C_2^d$ acts freely on $T^n$ if and only if there is $1$ in the sum of any distinct collection of rows of $E$. In this case $C_2^d$ is the holonomy group of the flat manifold $T^n/C_2^d$ if and only if there is either $2$ or $3$ in the sum of any distinct collection of rows of $E$ (see \cite[Lemma 2.4]{LPPS19}).

Let us consider the linear forms $\alpha,\beta\colon P\to\mathbb \F_2$ given by the following table

\begin{center}
    \begin{tabular}{|c|c|c|c|c|} \hline
        &0&1&2&3 \\ \hline\hline
        $\alpha$&0&1&1&0\\ \hline
        $\beta$&0&1&0&1\\ \hline
    \end{tabular}
\end{center}

Let $C_2^d\subseteq {D}^n$ and $j \in \{1,2,\ldots,n\}$. We define epimorphisms
\begin{equation}
\alpha_j:C_2^d\subseteq{D}^n\xrightarrow{pr'_j} {P}\xrightarrow{\alpha}\F_2,\;\;
\beta_j:C_2^d\subseteq{D}^n\xrightarrow{pr'_j} {P}\xrightarrow{\beta}\F_2,
\end{equation}
where
\[
pr'_j(g_{i_1},g_{i_2},\ldots,g_{i_n}) = i_j
\]
for $g_{i_1},g_{i_2},\ldots,g_{i_n} \in D$.
Since $H^1(C_2^d,\F_2)=\operatorname{Hom}(C_2^d,\mathbb \F_2)$, we can view $\alpha_j$ and $\beta_j$ as cohomology classes and define
\begin{equation}
\label{wzor33}
\theta_j=\alpha_j\cup \beta_j\in H^2(C_2^d,\F_2),\end{equation}
where $\cup$ denotes the cup product. Fix generators $b_1,\ldots,b_d$ of $C_2^d$. It is well known that
$H^*(C_2^d,\F_2)\cong \F_2[x_1,\ldots,x_d]$ where $\{x_1,\ldots,x_d\}$ is a basis of $H^1(C_2^d,\F_2)$ such that $x_i(b_j) = \delta_{ij}$ is the Kronecker delta.
Hence, elements $\alpha_j$ and $\beta_j$ correspond to
\begin{equation}\label{wzor34}
\alpha_j=\sum_{i=1}^d\alpha(pr_j(b_i))x_i,\;\;\; \beta_j=\sum_{i=1}^d \beta(pr_j(b_i))x_i\in C_2^d[x_1,\ldots,x_d],
\end{equation}
(see \cite[{Theorem 1.2 and Proposition 1.3}]{CMR10}). Moreover, if $[p_{ij}]$ is a $P$-matrix of $C_2^d$, which corresponds to the chosen generators,
we can write equations (\ref{wzor33}) and (\ref{wzor34}) as follows
\begin{equation}\label{wzor35}
\alpha_j=\sum_{i=1}^d \alpha(p_{ij})x_i,\;\; \beta_j=\sum_{i=1}^d \beta(p_{ij})x_i,\;\;\theta_j=\alpha_j\beta_j.
\end{equation}

There is an exact sequence
$$
0\to H^1(C_2^d,\F_2)\xrightarrow{\pi^*}H^1(\Gamma,\F_2)\xrightarrow{\iota^*}H^1(\Z^n,\F_2)\xrightarrow{d_2}
H^2(C_2^d,\F_2)\xrightarrow{\pi^*}H^2(\Gamma,\F_2)$$
where $d_2$ is the transgression and $\pi^*$ is induced by the quotient map $\pi:\Gamma\to C_2^d$ (see \cite{E91}).
\begin{prop}[{\cite[Proposition 3.2]{LPPS19}}]
    \label{Prop32}
    Suppose that free and diagonal action of $C_2^d$ on $T^n$ corresponds to a $P$-matrix $E$, which defines elements $\alpha_j,\beta_j$ and $\theta_j$ as in \eqref{wzor35}, for $1 \leq j \leq n$. Let $M=T^n/C_2^d$ and $\Gamma=\pi_1(M)$ be associated to the group extension (\ref{wzor21}). Then
    \begin{enumerate}
        \item $\theta_l=d_2(\varepsilon_l)$, where $\{\varepsilon_1,\ldots,\varepsilon_n\}$ is the basis of $H^1(\Z^n,\F_2)$ dual to the standard basis of $\Z^n \otimes_{\Z} \F_2$, for $l=1,\ldots,n$;
        \item the total Stiefel-Whitney class of $M$ is given by
        $$w(M)=\pi^*\left(w\right)\in H^*(\Gamma,\F_2)=H^*(M,\F_2),$$
        where
    \begin{equation}\label{wzor37}
    w=\prod_{j=1}^n\left(1+\alpha_j+\beta_j\right)\in\F_2[x_1,\ldots,x_d].
    \end{equation}
    \end{enumerate}
\end{prop}

We call the ideal
$$\Id_E = \langle\operatorname{Im}(d_2)\rangle=\langle \theta_1,\ldots,\theta_n\rangle\subseteq\F_2[x_1,x_2,\ldots,x_n]$$
the \emph{characteristic ideal} of $E$ and the quotient $C_E=\F_2[x_1,\ldots,x_d]/\Id_E$ -- the \emph{characteristic algebra} of $E$.

\begin{cor}[{\cite[Corollary 3.3]{LPPS19}}]
    \label{cor:total_sw_class}
    Suppose that free and diagonal action of $C_2^d$ on $T^n$ corresponds to a $P$-matrix $E$. There is a canonical homomorphism of graded algebras $\Phi:C_E\to H^*(T^n/C_2^d,\F_2)$ such that $\Phi([w])=w(T^n/C_2^d)$. Moreover, $\Phi$ is a monomorphism in degree less that or equal to two.
\end{cor}
\begin{defi}
    Given a $P$-matrix $E\in P^{d\times n}$, we define the Stiefel-Whitney class of $E$, to be the class $[w]\in C_E$ defined by (\ref{wzor37}).
\end{defi}
\begin{cor}
\label{cor:spin_structure}
Using notation of Corollary \ref{cor:total_sw_class}, let $w_2$ be the sum of degree $2$ summands of $w$. Then, by \cite[Proposition on page 40]{F00}, $T^n/C_2^d$ admits a spin structure if and only if $w_2 \in \Id_E$.
\end{cor}

Now, we describe a real Bott manifold $M(A)$ for strictly upper triangular matrix $A=[a_{ij}]$ with entries 0 or 1.
For every $1 \leq i \leq n$ we define an element $s_i$ of the group $\operatorname{Iso}(\R^n) = O(n) \ltimes \R^n$ of isometries of the euclidean space $\R^n$ in the following way. If $i \neq n$ then
\begin{equation}\label{gener}
s_i=\left(\operatorname{diag}\left[1,\ldots,1,(-1)^{a_{i,i+1}},\ldots,(-1)^{a_{i,n}}\right],\left(0,\ldots,0,\frac12,0\ldots,0\right)^T\right)
\end{equation}
where $(-1)^{a_{i,i+1}}$ is at the $(i+1, i+1)$ position and $\frac{1}{2}$ is on the $i$-th coordinate of the translational part.
In addition $s_{n} = \left(I,\left(0,0,\ldots,0,\frac12\right)\right)$. The group $\Gamma(A)$ generated by $s_1,\ldots,s_n$ is a crystallographic group. The subgroup generated by
$s_{1}^{2},s_{2}^{2},...,s_{n}^{2}$ consists of all translations of $\Z^n$. The action of $\Gamma(A)$ on $\R^n$ is free and the orbit space $M(A) = \R^n/\Gamma(A)$ is a flat manifold -- a closed Riemannian manifold with zero sectional curvature.

\section{Main results}
\label{sec:main}

We keep the notation of the previous section. Using the same methods as in \cite{PS16} and \cite{LPPS}, for each strictly upper triangular matrix $A=[a_{ij}]$ which generates the fundamental group of real Bott manifold $M(A)$ we get that the corresponding $P$-matrix $P_A=[p_{ij}]$ is of the form
\begin{equation}
\label{eq:pa}
P_A = \begin{bmatrix}
1 & 2a_{12} & \ldots & 2a_{1,n-1} & 2a_{1n}\\
0 & 1       & \ldots & 2a_{2,n-1} & 2a_{2n}\\
      & & \ddots & \\
0 & 0 & \ldots & 1 & 2a_{n-1,n}\\
0 & 0 & \ldots & 0 & 1
\end{bmatrix}.
\end{equation}
Note that in the above notation $2a_{ij}$ is multiplication in integers. To be more specific, we have
\[
p_{ij} = \left\{
\begin{array}{ll}
    1 & \text{ if } i=j,\\
    0 & \text{ if } i \neq j \text{ and } a_{ij}=0,\\
    2 & \text{ if } i \neq j \text{ and } a_{ij}=1,\\
\end{array}
\right.
\]
for $1 \leq i,j \leq n$.
Using the form \eqref{eq:pa} of $P_A$ and definition of forms $\alpha$ and $\beta$ we get that
\[
\begin{aligned}
\alpha_j &=\sum_{i=1}^{j}\alpha\left(p_{ij}\right)x_i
=\sum_{\substack{i=1}}^{j-1}\alpha(p_{ij})x_i+\alpha(p_{jj})x_j=\sum_{\substack{i=1}}^{j-1}a_{ij}x_i+x_j,\\
\beta_j & =\sum_{i=1}^j\beta\left(p_{ij}\right)x_i=\beta\left(p_{jj}\right)x_j=x_j
\end{aligned}
\]
and hence
\begin{align}
\alpha_j+\beta_j&=\sum_{i=1}^{j-1}a_{ij}x_i, \label{eq:a+b}\\
\theta_j &= \sum_{\substack{i=1}}^{j-1}a_{ij}x_ix_j+x_j^2 = \sum_{i \neq j}a_{ij}x_ix_j+x_j^2\label{eq:theta}
\end{align}
for all $1\leq j \leq n$. Note that in the last formula we take advantage of the definition of a Bott matrix.

We go back to the matrix $A$ of the real Bott manifold $M(A)$. From
{\cite{I11}} we have the following necessary and
sufficient condition for the existence of a K\"{a}hler structure on $M(A)$
\begin{theo}[{\cite[Theorem 3.1]{I11}}]\label{thI}
 Let $A$ be $2n$-dimensional matrix of real Bott manifold $M(A)$. Then the following conditions are equivalent:
   \begin{enumerate}
     \item there exist $n$ subsets $\{j_1,j_{n+1}\} \ldots \{j_n,j_{2n}\}$ of the set $\{1,2,\ldots,2n\}$ such that
       \[\bigcup_{k=1}^n\{j_k,j_{k+n}\}=\{1,2,\ldots,2n\}\]
       and $A^{(j_k)}=A^{(j_{k+n})}$ for all $1\leq k \leq n$, where $A^{(i)}$ is the $i$-th column of the matrix $A$;
       \item there exist a K\"{a}hler structure on $M(A)$.
  \end{enumerate}
\end{theo}

{Let $A$ be a Bott matrix of an RBK-manifold.}
Using the notation from Theorem {\ref{thI}}, let $\tilde A=[A^{(j_1)}A^{(j_2)}\ldots A^{({j_n})}]$ be a matrix obtained from $A$ {by removing duplicated columns}. Then $\tilde A \in \F_2^{2n\times n}$. Let
$$\tilde{S}_i=\sum_{k=1}^na_{ij_k} \in \F_2$$
denote the sum of elements in the $i$-th row of the matrix $\tilde A$, for $1\leq i\leq 2n$. Let us recall our main theorem.

\begin{mtheo}\label{theo1}
    Let $A$ be a matrix of of $2n$-dimensional RBK-manifold $M(A)$. Then $M(A)$ admits a spin structure if and only if for every $1 \leq i \leq 2n$
$$
\tilde{S}_i=1\Longrightarrow A^{(i)}=0.
$$
\end{mtheo}

\noindent{\rm Proof.} $M(A)$ is $2n$-dimensional RBK-manifold, so from \eqref{eq:a+b}
we get
$$\alpha_{j_{k+n}}+\beta_{j_{k+n}}=\alpha_{j_k}+\beta_{j_k}=\sum_{i=1}^{j_k-1}a_{ij_{k}}x_i,$$
and
$$\begin{aligned}
w(M(A))&=\prod_{k=1}^{2n}\left(1+\alpha_{j_k}+\beta_{j_k}\right)=\prod_{k=1}^{n}\left(1+\alpha_{j_k}+\beta_{j_k}\right)^2\\
&=\prod_{k=1}^{n}\left(1+\sum_{i=1}^{j_k-1}a_{ij_k}x_i\right)^2=\prod_{k=1}^{n}\left(1+\sum_{i=1}^{j_k-1}a_{ij_k}x_i^2\right).\\
\end{aligned}$$
From the above considerations and from the definition of a Bott matrix we have
$$\begin{aligned}
w_2 = w_2(M(A))&=
\sum_{k=1}^n\sum_{i=1}^{j_k-1}a_{ij_k}x_i^2=\sum_{k=1}^n\sum_{i=1}^{2n}a_{ij_k}x_i^2\\
&=\sum_{i=1}^{2n}\sum_{k=1}^n a_{ij_k}x_i^2=\sum_{i=1}^{2n}\tilde{S}_ix_i^2.
\end{aligned}$$
Let $J=\{j:\tilde{S}_j=1\} \subset \{1,\ldots,2n\}$. Then
\begin{equation}\label{wzor111}
w_2=\sum_{j\in J}x_j^2.
\end{equation}
By Corollary \ref{cor:spin_structure}, the existence of a spin structure on $M(A)$ is equivalent to $w_2 \in \Id_{P_A}$, which -- by \eqref{eq:theta} -- occurs if and only if
\begin{equation}\label{wzor222}
\sum_{j\in J}x_j^2=\sum_{j\in J}\theta_j.
\end{equation}
The above holds if $A^{(j)}=0$, since then, using formula \eqref{eq:theta} again, we get $\theta_j = x_j^2$, for $j \in J$.

Now, let
$$
K=\{(i,j):a_{ij} =1, j\in J, 1 \leq i \leq 2n\}.
$$
Then
\begin{equation}
\label{eq:sum_of_thetas}
\sum_{j\in J}\theta_j=\sum_{j\in J}\left(\sum_{i\neq j}a_{ij}x_ix_{j}+x_j^2\right)=\sum_{j\in J}x_j^2+\sum_{(i,j)\in K}x_ix_{j}.
\end{equation}

From \eqref{wzor111}, \eqref{wzor222} and \eqref{eq:sum_of_thetas} we have
$$\sum_{j\in J}x_j^2=\sum_{j\in J}x_j^2+\sum_{(i,j)\in K}x_ix_{j},$$
hence $K=\emptyset$, which means that the $j$-th column of the Bott matrix $A$ has only zero entries.

\hskip13cm$\square$

\begin{ex}
    Let $$A=\left[\begin{matrix}0&0&1&1&1&1\\0&0&1&1&1&1\\0&0&0&0&1&1\\0&0&0&0&1&1\\0&0&0&0&0&0\\0&0&0&0&0&0\end{matrix}\right]$$
    be a matrix of a manifold $M(A)$. Then
$$
\tilde A=\left[\begin{matrix}0&1&1\\0&1&1\\0&0&1\\0&0&1\\0&0&0\\0&0&0\end{matrix}\right],
$$
$\tilde{S}_1=\tilde{S}_2=\tilde{S}_5=\tilde{S}_6=0$, $\tilde{S}_3=\tilde{S}_4=1$ and there are entries equal to $1$ in columns $A^{(3)}$ and $A^{(4)}$, so $M(A)$ has no spin-structure.
\end{ex}

At the end let us note an easy corollary of our main theorem.

\begin{cor}
Let $A \in \F_2^{2n \times 2n}$ be a Bott matrix such that $M(A)$ is a RBK-manifold. Let $T = \{ j : A^{(j)} \neq 0 \}$. If
\[
T = T_1 \sqcup \ldots \sqcup T_l
\]
such that for every $1 \leq k \leq l$ $T_k$ is a four-elements set and
\[
\forall_{i,j \in T_k} A^{(i)}=A^{(j)},
\]
then $M(A)$ admits a spin structure.
\end{cor}

\bibliographystyle{plain}
\bibliography{bibl}

\end{document}